\documentclass[10pt, oneside, reqno]{amsart}
\usepackage[utf8]{inputenc}
\usepackage[english]{babel}
\usepackage[T1]{fontenc}
\usepackage{amsmath}
\usepackage{amsfonts}
\usepackage{amssymb}
\usepackage{amsthm}
\usepackage{amscd}
\usepackage{soul}
\usepackage{geometry}
\usepackage{enumitem} 
\usepackage{cancel} 
\usepackage{graphicx}
\usepackage{mathtools}
\usepackage{color}
\usepackage{comment}
\usepackage{mathdots}
\usepackage{algorithm}
\usepackage{algpseudocode}
\usepackage{tikz}
\usetikzlibrary{calc,arrows.meta,decorations.pathreplacing}
\usetikzlibrary{arrows}
\usepackage{enumitem}
\usepackage{color}
\usepackage{bbm}
\definecolor{marin}{rgb}   {0.,   0.3,   0.7} 
\definecolor{rouge}{rgb}   {0.8,   0.,   0.} 
\definecolor{sepia}{rgb}   {0.8,   0.5,   0.} 
\definecolor{uuuuuu}{rgb}{0.26666666666666666,0.26666666666666666,0.26666666666666666}
\usepackage[colorlinks,citecolor=marin,linkcolor=rouge,
            bookmarksopen,
            bookmarksnumbered
           ]{hyperref}

\usetikzlibrary{calc,backgrounds,fit}
\usetikzlibrary{external}
\usetikzlibrary{patterns}
\usetikzlibrary{decorations.pathreplacing,decorations.markings,arrows}
\usetikzlibrary{shapes,intersections,positioning}
\usetikzlibrary{arrows}

\newcommand\N{\mathbb{N}}

\newcommand\R{\mathbb{R}}
\newcommand\C{\mathbb{C}}

\newcommand{\dd}{\mathrm{d}}

\newcommand{\enstq}[2]{\left\{#1~\middle|~#2\right\}}

\newcommand\eps{\varepsilon}

\renewcommand{\Im}{\operatorname{Im}}

\newcommand\Id{\mathrm{Id}}

\newcommand\exterior{\mathrm{ext}}
\newcommand\interior{\mathrm{int}}
\newcommand\diam{\mathrm{diam}}

\newcommand\err{\mathrm{err}}
\newcommand\centre[4][below]{\node (#3) at #2 [circle,minimum size=0.5em,inner sep=0pt,thin,fill,solid] {}; \node [#1=0.1em] at (#3) {#4};}

\newtheorem{theorem}{Theorem}

\newtheorem{corollary}{Corollary}
\newtheorem{lemma}{Lemma}
\newtheorem{definition}{Definition}

\title[Finite Volumes for the Gross-Pitaevskii equation]{Finite Volumes for the Gross-Pitaevskii equation}

\author{Quentin Chauleur}
\address{INRIA Lille, Univ Lille \& Laboratoire Paul Painlevé,
CNRS UMR 8524 Lille, Cité Scientifique, 59655 Villeneuve-d'Ascq, France. }
\email{Quentin.Chauleur@math.cnrs.fr}


\begin{document}

\maketitle

\begin{abstract}
We study the approximation by a semi-discrete finite-volume scheme of the Gross-Pitaevskii equation with time-dependent potential in two dimensions, performing a two-point flux approximation scheme in space. We rigorously analyze the error bounds relying on discrete uniform Sobolev inequalities. We finally perform some numerical simulations to investigate convergence error.
\end{abstract}

\section{Introduction}
Let $\Omega$ be an open bounded convex polygonal subset of $\R^2$ with Lipchitz boundary $\partial \Omega$. For a time horizon $T>0$, we consider the Gross-Pitaevskii equation 
\begin{equation} \label{gross_pitaevskii} \tag{GP}
i  \partial_t u + \Delta u= \gamma |u|^2 u + V(t,x) u
\end{equation}    
in $(0,T) \times \Omega$ with homogeneous Dirichlet boundary conditions $u(t,x)=0$ for all $(t,x)\in \left[0,T\right] \times \partial \Omega$, defocusing nonlinearity $\gamma>0$ and initial condition $u(0)=u_0$. Equation \eqref{gross_pitaevskii} appears as a fundamental model to describe the evolution of Bose-Einstein condensates, a state of matter near absolute zero temperature whose dynamics can be described by a single wave function $u$. The time-dependent potential~$V=V(t,x)\in \R$ then both act as a magnetic confining potential and a stirring perturbation induced by laser beams. When set to rotation in a strong confinement regime, Bose-Einstein condensates exhibit complex nonlinear phenomena such as quantum turbulence, which have drawn more and more attention from the scientific community over the last decades. In particular, the nucleation of vortices with a quantized circulation, which are related to superfluid properties, as well as their interactions, played an important role for irreversible energy transfer mechanisms in quantum turbulent fluids~\cite{krstulovic2020}.

The recent developments of laser trapping and cooling has made BEC experiments very precise in a number of complex configurations with various geometries. On the other hand, numerical simulations are usually made on periodic fields to allow the use of the Fast Fourier Transform in order to retain spectral accuracy, or with finite differences \cite{bao2013} on square or cubic lattices. This motivates the development of new methods in order to efficiently simulate the Gross-Pitaevskii equation on more general geometries. Here we will only assume that $\Omega$ is an open bounded convex polygonal set, and we analyze a numerical scheme based on a finite volume approximation in space. We perform a Two-Point Flux Approximation (TFPA) finite volume scheme, which stands as a very popular method in numerous applications, as it is both straightforward to implement and robust. Note that to the best of the author's knowledge, this work seems to be one of the first contribution to the numerical analysis of a finite volume scheme in the context of nonlinear dispersive PDEs. Even in the case of linear Schrödinger equations, the only work we are aware of is the paper of Bradji \cite{bradji2015}, in the context of hybrid finite volume scheme on general nonconforming meshes.

Finite elements methods has also recently been developed for the simulation of quantum fluids, in particular with the works of Henning and M\r{a}lqvist \cite{henning2017} on convex domains, with GPElab developed by Antoine and Duboscq \cite{duboscq2015stationary,duboscq2015dynamics} or the toolbox developed by Vergez, Danaila, Auliac and Hecht \cite{danaila2016} on the FreeFEM++ software. Finite element methods are based on variational formulas, and usually require an adequate functional framework that we avoid here with our finite volume approach.

From the theoretical point of view, the only work dealing with the Cauchy problem of nonlinear Schrödinger with time dependent potential seems to be \cite{carles2011}. In particular this paper ensures the existence of a unique global solution of equation \eqref{gross_pitaevskii} on the whole space $\R^d$ for $d \leq 3$, under the condition that $V(t,x)$ is real-valued, locally bounded in time and subquadratic in space. Also note that equation \eqref{gross_pitaevskii} satisfies several conservation laws, such as the mass conservation
\[ \|u(t)\|_{L^2(\Omega)}=\|u_0\|_{L^2(\Omega)}   \]
for all $t \in \R$, or the energy balance law
\[ \frac{\dd }{\dd t} H(t) = \int_{\Omega} \partial_t V(t,x) |u(t,x)|^2 \dd x   \]
where
\[ H(t)\coloneqq \| \nabla u \|_{L^2(\Omega)}^2 + \frac{\gamma}{2} \| u(t) \|_{L^4(\Omega)}^4 +  \int_{\Omega}V(t,x) |u(t,x)|^2 \dd x. \]
In particular, if $V=0$, the energy $H$ is also conserved over time. Throughout all this paper we will make the assumption on the potential that $V \in \mathcal{C}^1(\left[ 0, T \right], H^3(\Omega))$, that $V$ is bounded from below, and that $\partial_t V$ is bounded from above. We also suppose that these assumptions would reasonably induce the global existence of a unique smooth solution $u \in \mathcal{C}(\left[ 0, T \right], H^5\cap H^1_0(\Omega))$ with initial condition $u_0 \in H^5\cap H^1_0(\Omega)$ on the open bounded subset $\Omega$ of $\R^d$.

This paper is organized as follows. In section \ref{section_main_result}, we will recall some notations about finite volume schemes in order to state our main result Theorem \ref{theorem_convergence}. Section \ref{section_discrete_functional_analysis} introduces some useful discrete functional properties along discrete conserved quantities of the numerical scheme. In section \ref{section_error}, we show several error estimates which allows to prove the convergence result from Theorem \ref{theorem_convergence}. Finally, in section \ref{section_numerics} we perform some simulations to investigate numerical convergence.

\section{Discrete framework and main result} \label{section_main_result}

\subsection{Finite volume notations}
We consider finite-volumes admissible meshes in the sense of the following definition, as introduced in the seminal book of Eymard, Gallouët and Herbin \cite{eymard2000}.

\begin{definition} \textbf{(Admissible finite volume discretization)}. \label{defmesh}  \\
An admissible finite-volume discretization of $\Omega$ is a pair $\mathcal{V}=(\mathcal{M},\mathcal{E})$ such that:
\begin{itemize}
\item The mesh $\mathcal{M}$ (see Figure~\ref{fig:notation_mesh}) is given by a family of open polygonal and convex subsets $K$ of~$\Omega$, called \textit{control volumes} of $\mathcal{T}$, satisfying the following properties:
\begin{itemize}
\item $\overline{\Omega}=\bigcup_{K\in\mathcal{M}}\overline{K}$.
\item If $K,L\in\mathcal{M}$ with  $K\neq L$ then $\operatorname{int}K\cap\operatorname{int}L=\emptyset$.
\item If $K,L\in\mathcal{M}$, with $K\neq L$ then either the $1$-dimensional Lebesgue measure of $\overline{K}\cap \overline{L}$ is~$0$ or $\overline{K}\cap \overline{L}$ is the edge of the mesh denoted $\sigma=K|L$  separating the control volumes~$K$ and $L$.
\item To each control volume $K\in\mathcal{M}$, we associate a point $x_K\in \overline{K}$ (called the center of $K$) such that if $K,L\in\mathcal{M}$ are two neighboring control volumes, the straight line between the centers $x_K$ and $x_L$ is orthogonal to $\sigma=K|L$. 
\end{itemize}
\item The set of edges $\mathcal{E}$ is a partition of the mesh skeleton $\cup_{K \in \mathcal{M}} \partial K$ into edges $\sigma$ which are subsets contained in hyperplanes of $\overline{\Omega}$, such that $\partial \mathcal{M} = \cup_{\sigma \in \mathcal{E}} \overline{\sigma}$. We respectively define
\[ \mathcal{E}_{\interior} \coloneqq \enstq{\sigma \in \mathcal{E}}{\sigma \notin \partial \Omega} \quad \text{and} \quad \mathcal{E}_{\exterior} \coloneqq \enstq{\sigma \in \mathcal{E}}{\sigma \in \partial \Omega}  \]
the subsets of interior and exterior edges of $\mathcal{E}$. We also denote by $\mathcal{E}_K$ the set of edges of any $K \in \mathcal{M}$, as well as $\nu_{K,\sigma}\in \R^2$ the unit normal vector to $\sigma$ pointing outwards from $K$. We also denote by $\mathcal{M}_{\interior} \subset \mathcal{M}$ the set such that for all $K \in \mathcal{M}_{\interior}$ and all $\sigma \in \mathcal{E}_K$, $\sigma \notin \partial \Omega$, and we define $\mathcal{M}_{\exterior}=\mathcal{M} \backslash \mathcal{M}_{\interior}$.
\end{itemize}
We denote by $|K|\coloneqq\lambda_2(K)$ the area of the control volume $K$, where $\lambda_2$ stands for the $2$-dimensional Lebesgue measure. We also denote by $h_K=\diam(K)=\sup \enstq{|x-y|}{x,y \in K}$ the diameters of control volumes $K \in \mathcal{M}$, and by
\[ h_{\mathcal{M}} \coloneqq \sup_{K \in \mathcal{M}} h_K \]
the mesh size of $\mathcal{V}$. Moreover, for $K$ and $L$ two neighboring control volumes of $\mathcal{M}$, we denote by $\sigma = K|L \in \mathcal{E}_{\interior}$ their common edge, and by $|\sigma|=\lambda_{1}(\sigma)$ the length of $\sigma$. We also denote by $d_{K|L}$ the distance between $x_K$ and $x_L$, and by $d_{K,\sigma}$ the length of the segment between $x_K$ and $\sigma$, orthogonal to $\sigma$, if $\sigma \in \mathcal{E}_{\exterior}$. Finally, we denote by $d_{\mathcal{M}} \in \N$ the number of control volumes $K \in \mathcal{M}$ of the finite-volume mesh, and by $e_{\mathcal{E}}$ the number of edges of $\mathcal{V}$. Once again, we refer to Figure \ref{fig:notation_mesh} for an illustration of these notations.
\end{definition}

\begin{figure}[htbp!]
\centering
\begin{tikzpicture}[scale=2]

  \clip (-1.2,-0.6) rectangle (1.8,1.3);

  \node[rectangle,fill] (a) at (-1,0.6) {};
  \node[rectangle,fill] (b) at (0,1.2) {};
  \node[rectangle,fill] (c) at (0,-0.2) {};
  \node[rectangle,fill] (d) at (1.5,0.3) {};

  \centre[above right]{(-0.6,0.5)}{xK}{$x_K$};
  \centre[above left]{(0.9,0.5)}{xL}{$x_L$};
  
  \draw[thick] (b)--(c) node [pos=0.7,right] {$\sigma=$\small{$K|L$}};

  \draw[thin,opacity=0.5] (a) -- (b) -- (c) -- (a) ;
  \draw[thin,opacity=0.5] (d) -- (b) -- (c) -- (d);

  \draw[dashed] (xK) -- (xL);
  
  \coordinate (KK) at ($(xK)!0.65!-90:(xL)$);
  \coordinate (LL) at ($(xL)!0.65!90:(xK)$);

  \draw[dotted,thin] (xK) -- (KK);
  \draw[dotted,thin] (xL) -- (LL);

  \draw[|<->|] (KK) -- (LL) node [midway,fill=white,sloped] {$d_{K|L}$};
 
  \coordinate (KAB) at ($(a)!(xK)!(b)$);
  \coordinate (KAC) at ($(a)!(xK)!(c)$);

  \coordinate (LDB) at ($(d)!(xL)!(b)$);
  \coordinate (LDC) at ($(d)!(xL)!(c)$);

  \draw[dashed] (xK) -- ($(xK)!3!(KAB)$);
  \draw[dashed] (xK) -- ($(xK)!3!(KAC)$);

  \draw[dashed] (xL) -- ($(xL)!3!(LDB)$);
  \draw[dashed] (xL) -- ($(xL)!3!(LDC)$);

  \begin{scope}[on background layer]   
    \draw (0,0.5) rectangle ++ (0.1,-0.1);
  \end{scope}
  
  \coordinate (nukl) at ($(b)!0.3!(c)$);
  \draw[->,>=latex] (nukl) -- ($(nukl)!0.3cm!90:(c)$) node[above] {$\nu_{K,\sigma}$};

\end{tikzpicture}
\caption{Notations for the mesh $\mathcal{M}$ associated with $\Omega$\label{fig:notation_mesh}.}
\end{figure}

In the particular case of \textit{Voronoï meshes}, each control volume $K\in \mathcal{M}$ associated to a center $x_K$ is defined by
\[  K= \enstq{x \in \Omega}{ |x-x_K|<|x-x_L|, \ \forall x_L \in \mathcal{M}, \ x_L \neq x_K}.  \] 
In fact, regular Voronoï meshes often satisfy
\begin{equation} \label{eq_boundary_voronoi} \tag{boundary}
\mathcal{E}_K \cap \mathcal{E}_{\exterior}  \neq \emptyset \Rightarrow x_K \in \partial \Omega,
\end{equation} 
a property that we will always assume in the following. The use of the Voronoï method to produce admissible finite volume meshes goes back to the work of Mishev \cite{mishev1998} and has now proven being useful in a number of applications \cite{gartner2019}. As pointed in \cite[Example 9.2]{eymard2000}, an advantage of the Voronoï method is that it easily leads to meshes on non polygonal domain $\Omega$, which might be of interest for physical applications. Note also that a relationship can be established between TPFA scheme on Voronoï meshes and generalized mixed-hybrid mimetic finite difference method, and super-convergence can occur \cite{droniou2018}. 

In two dimension, a usual way to construct meshes satisfying the orthogonality property $(x_K,x_L) \perp \sigma$ is to partition $\Omega$ into a conforming triangulation with acute triangles called Delaunay triangulation, and to take each $x_K$ as the circumcenter of $K$. This method is one of the most popular in computational mesh generation, and the dual mesh then corresponds to a Voronoï admissible finite volume mesh. In our upcoming analysis, we will assume the following regularity constraint on the discretization $\mathcal{V}$: there exists $0<\Theta_{\mathcal{V}} \leq 1$ independent of $h_{\mathcal{M}}$ such that
\begin{equation} \label{eq_regularity_mesh} \tag{reg}
 h_K \Theta_{\mathcal{V}} \leq d_{K,\sigma} \leq h_K \quad \text{and} \quad \left\{ \begin{aligned}
& \Theta_{\mathcal{V}} h_K \leq |\overline{\sigma^T}| \leq h_K, \\
& \Theta_{\mathcal{V}} h_K \leq |\overline{\sigma^{T'}}| \leq h_K, 
\end{aligned} \right. 
\end{equation}
for any $K \in \mathcal{M}$ and $\sigma \in \mathcal{E}_K$, writing $\overline{\sigma}=\overline{\sigma^T} \cup \overline{\sigma^{T'}}$ (see Figure \ref{fig:notations_dual_triangulation} in section \ref{section_discrete_functional_analysis} for precise definition).

\subsection{Discrete unknowns}
For $k \in \N$ and a subset $X \subset \overline{\Omega}$, we denote by $\mathbb{P}^k$ the space of polynomial functions $X \rightarrow \C$ of degree at most $k$. We then introduce the discrete space  
\[ X_{\mathcal{V},0} \coloneqq \enstq{v_{\mathcal{M}}:\Omega \rightarrow \C}{\forall K \in \mathcal{M}, \left. v_{\mathcal{M}}\right|_K \in \mathbb{P}^0(K)}.  \]
For convenience, such functions may be identified with vectors fields writing that $v_{\mathcal{M}}=(v_K)_{K \in \mathcal{M}}$, where $\forall K \in \mathcal{M}$, $\forall x \in K$, $v_{\mathcal{M}}(x)=v_K$. We also introduce the following discrete norms on our Voronoï mesh embedded with Dirichlet boundary conditions.

\begin{definition} \textbf{(Discrete spaces)}.  \\
Let $1 \leq p < \infty$. We define the discrete Lebesgue norms 
\[ \| v_{\mathcal{M}} \|_{0,p}^p \coloneqq \sum_{K \in \mathcal{M}} |K| |v_K|^p  \quad \text{and} \quad  \| v_{\mathcal{M}} \|_{0,\infty} \coloneqq \sup_{K \in \mathcal{M}} |v_K|, \]
as well as the discrete semi-norm 
\[ |v_{\mathcal{M}}|_{1,2,\mathcal{M}}^2 \coloneqq \sum_{K\in \mathcal{M}} \sum_{\sigma \in \mathcal{E}_K \cap \mathcal{E}_{\interior}} |\sigma| d_{K,\sigma} \left| \frac{D_{\sigma} v_{\mathcal{M}}}{d_{\sigma}} \right|^2 = \sum_{\sigma \in \mathcal{E}_{\interior}} \frac{|\sigma|}{d_{K|L}} |v_K-v_L|^2  \]
for all $v_{\mathcal{M}} \in X_{\mathcal{V},0}$, where
\[ D_{\sigma} v_{\mathcal{M}} = |v_K-v_L|  \quad \text{for} \ \sigma=K|L \in \mathcal{E}^{\interior}  \]
and $d_{\sigma}=d_{K,\sigma} + d_{L,\sigma}$ for $\sigma=K|L \in \mathcal{E}^{\interior}$.
\end{definition}

Note that there is no contribution from $\sigma \in \mathcal{E}_{\exterior}$ in the $|\cdot|_{1,2,\mathcal{M}}$ semi-norm as we specifically work on Voronoï meshes with homogeneous Dirichlet boundary conditions, which imposes $v_K=0$ if $K \in \mathcal{M}_{\exterior}$. Also, from the \textit{discrete Poincaré inequality} (see for instance \cite[Lemma 9.1]{eymard2000})
\[  \| v_{\mathcal{M}} \|_{0,2} \leq \diam(\Omega) | v_{\mathcal{M}} |_{1,2,\mathcal{M}}, \]
so that the semi-norm $|\cdot |_{1,2,\mathcal{M}}$ is in fact a norm. Finally we define the \textit{pointwise evaluation operator}~$P_{\mathcal{M}}:\mathcal{C}(\Omega) \rightarrow X_{\mathcal{V},0}$ of any continuous function $f \in \mathcal{C}(\Omega)$ by
\[ \left. P_{\mathcal{M}} f\right|_K  = f(x_K), \quad \forall K \in \mathcal{M}. \]

\subsection{Finite volume scheme and main result}
We denote by $A_{\mathcal{V}}:X_{\mathcal{V},0} \rightarrow X_{\mathcal{V},0}$ the discrete Laplace operator defined for all $K\in \mathcal{M}$ and $v_{\mathcal{M} }\in X_{\mathcal{V},0}$ by
\begin{equation*} \label{definition_discrete_laplacian_Dirichlet}
\left. A_{\mathcal{V}} v_{\mathcal{M} }\right|_K \coloneqq \left\{
\begin{aligned}
 & \frac{1}{|K|}\sum_{\sigma \in \mathcal{E}_K} \frac{|\sigma|}{d_{K|L}} (v_L-v_K) & \quad \text{if} \ K \in \mathcal{M}_{\interior}, \\
 & 0 & \quad \text{if} \ K \in \mathcal{M}_{\exterior}. 
 \end{aligned} \right.
\end{equation*}
We then define our numerical scheme $t \mapsto u_{\mathcal{M}}(t)$ for all $t\in \left[0,T\right]$ as the unique solution of the system
\begin{equation} \label{FVGP} \tag{FVGP}
\left\{
\begin{aligned}
& i u_{\mathcal{M}}' + A_{\mathcal{V}}u_{\mathcal{M}}=\gamma |u_{\mathcal{M}}|^2 u_{\mathcal{M}} + (P_{\mathcal{M}} V) u_{\mathcal{M}}, \\
&   u_{\mathcal{M}}(0)=P_{\mathcal{M}} u_0, 
\end{aligned} \right. 
\end{equation}  
Note that from the definition of $A_{\mathcal{V}}$ and $P_{\mathcal{M}}$, and as $u_0(x)=0$ for all $x \in \partial \Omega$, we directly get that~$u_{K}(t)=0$ for all $t \in \left[0,T \right]$ if $K \in \mathcal{M}_{\exterior}$. We now state the main result of this work:

\begin{theorem} \label{theorem_convergence}
Let $\mathcal{V}=(\mathcal{M},\mathcal{E})$ be an admissible finite-volume Voronoï tessellation of $\Omega$ satisfying the regularity properties \eqref{eq_boundary_voronoi}-\eqref{eq_regularity_mesh}. Let $u_0 \in H^5\cap H^1_0(\Omega)$ and $T>0$. We assume that there exists a unique solution $u \in \mathcal{C}(\left[0,T\right];H^5\cap H^1_0(\Omega))$ of the Gross-Pitaevskii equation \eqref{gross_pitaevskii} with $u(0)=u_0$. Let $u_{\mathcal{M}}(0)=P_{\mathcal{M}} u_0\in X_{\mathcal{V},0}$, and we denote by $u_{\mathcal{M}}=(u_K)_{K\in\mathcal{M}}$ the numerical scheme defined by equation \eqref{FVGP}. Then for any $\eps>0$, we have
\[  \sup_{t\in\left[0,T\right]} \| P_{\mathcal{M}} u(t)-u_{\mathcal{M}}(t) \|_{0,2} \leq C h^{\frac12-\eps}   \]
where $C=C(T,\Omega,\|u\|_{L^{\infty}_T H^5(\Omega)},\|V\|_{L^{\infty}_T H^3(\Omega)},\Theta_{\mathcal{V}},\eps)>0$.
\end{theorem}

It should be pointed out that with few adaptations, our result extends to Neumann boundary conditions $\nabla u \cdot \overrightarrow{\nu}$ for all $(t,x)\in \left[0,T \right] \times \partial \Omega$. It should also apply on rectangular meshes in the sense of \cite{Boyer2006}, since the key estimate from Corollary \ref{discrete_infty_lemma} is also valid on such meshes (see  \cite[Remark 2.1]{Boyer2006}). For clarity, however, we have chosen to make our statement as simple as possible.

We end this section with some notations. The uniform norm on the time interval $\left[0, T \right]$ will be denoted $\|\cdot\|_{L^{\infty}_T}$, for conciseness purposes. If not mention otherwise, an edge $\sigma \in \mathcal{E}_{\interior}$ is systematically associated to control volumes $K$ and $L \in \mathcal{M}$ such that $\sigma = K|L$ in summations. We will use the notation $\lesssim$ for estimates which does not depend on the parameter $h_{\mathcal{M}}$ or the Lebesgue exponent $1\leq p < \infty$.

\section{Discrete functional analysis} \label{section_discrete_functional_analysis}

\subsection{Voronoï tessellation and dual Delaunay triangulation}
From now on,  we consider a Voronoï tessellation $\mathcal{V}=(\mathcal{M},\mathcal{E})$ of $\Omega$ satisfying the regularity properties \eqref{eq_boundary_voronoi}-\eqref{eq_regularity_mesh}. The corresponding dual Delaunay triangulation is denoted $\mathcal{D}=(\mathcal{T},\mathcal{F})$, and is obtained connecting the cell centers of $\mathcal{V}$.  The resulting triangulation satisfy $\overline{\Omega}=\cup_{T \in \mathcal{T}} \overline{T}$ and $\partial \mathcal{T}= \cup_{F \in \mathcal{F}} \overline{F}$. Each Delaunay triangle $T\in\mathcal{T}$  can be split into 6 right-angled triangle denoted $T_{mn}$ for $1 \leq m\neq n \leq 3$. For two neighboring triangles $T$, $T'\in\mathcal{T}$ such that $\overline{\sigma}=\left[x_T,x_{T'}\right]$, we define $\sigma^T=(x_T,x_{\sigma})$ and $\sigma^{T'}=(x_{T'},x_{\sigma})$ such that $\overline{\sigma}=\overline{\sigma^T} \cup \overline{\sigma^{T'}}$. We refer to Figure \ref{fig:notations_dual_triangulation} for these auxiliary notations.

\begin{figure}[htbp!]
\centering
\begin{tikzpicture}[scale=3,>=Latex]

\coordinate (K1) at (1.73,0);
\coordinate (K2) at (0,1);
\coordinate (K3) at (0,-1);
\coordinate (center) at (0.58,0);
\coordinate (sigma12) at (0.87,0.5);
\coordinate (xT) at (1.07,0.85);

\draw[thick] (K1)--(K2)--(K3)--cycle;

\draw[dashed] (K1)--(center) node[midway,below] {$T_{13}$};
\draw[dashed] (K2)--(center) node[midway,left] {$T_{23}$};
\draw[dashed] (K3)--(center) node[midway,left] {$T_{32}$};
\draw[dashed] (K1)--(sigma12) node[midway,below] {$T_{12}$};
\draw[dashed] (K2)--(sigma12) node[midway,below] {$T_{21}$};
\draw[dashed] (K3)--(sigma12) node[midway,below right] {$T_{31}$};

\fill[blue] (K1) circle (0.8pt) node[below right] {$x_{K_1}$};
\fill[blue] (K2) circle (0.8pt) node[above left] {$x_{K_2}$};
\fill[blue] (K3) circle (0.8pt) node[below left] {$x_{K_3}$};
\fill (center) circle (0.8pt) node[above left] {$x_T$};
\fill (sigma12) circle (0.8pt) node[right] {$x_{\sigma_{12}}$};
\fill (xT) circle (0.8pt) node[above,left] {$x_{T'}$};

\draw[blue] (center) -- ++(-1,0);
\draw[blue] (center) -- ++ (0.5,-0.83);
\draw[blue] (center) -- (xT);
\draw[blue] (xT) -- ++(-0.3,0.5);
\draw[blue] (xT) -- ++(0.6,0);

\fill[blue] (1.5,-0.5) node {$K_1$};
\fill[blue] (-0.5,0.5) node {$K_2$};
\fill[blue] (-0.5,-0.5) node {$K_3$};

\draw[<-] (0.77,0.3) -- (1.4,0.4) node[right] {$\sigma_{12}^{T} \coloneqq (x_T,x_{\sigma_{12}})$};
\draw[<-] (1,0.7) -- (1.4,0.65) node[right] {$\sigma_{12}^{T'}$};
\draw[<->] (0.01,-1.05) -- (0.89,-0.54) node[midway, below right] {$d_{K_3,\sigma_{13}}$};

\draw  +(-0.05,0.35) -- +(0.05,0.45);
\draw  +(-0.05,0.40) -- +(0.05,0.50);

\draw  +(-0.05,-0.45) -- +(0.05,-0.35);
\draw  +(-0.05,-0.50) -- +(0.05,-0.40);

\draw  +(-0.05,0.05) -- +(0,0.05);
\draw  +(-0.05,0.05) -- +(-0.05,0);

\end{tikzpicture}
\caption{Notations of the dual Delaunay triangulation $\mathcal{D}$.} \label{fig:notations_dual_triangulation}
\end{figure}

As a direct consequence of the mesh regularity assumptions, denoting $\mathcal{N} \coloneqq \enstq{(m,n)}{1 \leq m \neq n \leq 3}$, for all~$(m,n),(m',n') \in \mathcal{N}$, we have
\[ |T_{mn}| \lesssim |T_{m'n'}|.   \]
We introduce the space
\[ X_{\mathcal{D},0} \coloneqq \enstq{v_{\mathcal{T}} \in \mathcal{C}_0(\Omega)}{\forall T \in \mathcal{T}, \left. v_{\mathcal{T}}\right|_T \in \mathbb{P}^1(T)} \subset H^1_0(\Omega).  \]
To a discrete function $v_{\mathcal{M}} \in X_{\mathcal{V},0}$ on the primal Voronoï mesh we can associate a continuous function $v_{\mathcal{T}} \in X_{\mathcal{D},0}$ on the dual Delaunay mesh, which interpolates on each triangle $T \in \mathcal{T}$ with corners $x_{K_1}$, $x_{K_2}$ and $x_{K_3}$ the three values $v_{K_1}$, $v_{K_2}$ and $v_{K_3}$ of $v_{\mathcal{M}}$ on the associated Voronoï cells.

\begin{lemma} \label{lemma_eq_discrete_H1_dual}
There exists a constant $C=C(\Omega,\Theta_{\mathcal{V}})>0$ such that for all $v_{\mathcal{M}} \in X_{\mathcal{V},0}$,
\[ \| \nabla v_{\mathcal{T}} \|_{L^2(\Omega)} \leq C |v_{\mathcal{M}} |_{1,2,\mathcal{M}}.  \]
\end{lemma}
\begin{proof}
Let $T\in\mathcal{T}$. We use the notations from Figure \ref{fig:notations_dual_triangulation}. There is three possibilities:
\begin{enumerate}[label=\textnormal{(\roman*)}]
	\item $\overline{T}=\cup_{(m,n)\in \mathcal{N}} \overline{T_{m,n}}$, and $v_{\mathcal{T}}$ restricted on $T$ interpolates the values $v_{K_1}$, $v_{K_2}$ and $v_{K_3}$.
	\item $\overline{T}= \overline{T_{m,n}} \cup \overline{T_{n,m}}$ for some $(m,n)\in \mathcal{N}$, and $v_{\mathcal{T}}$ restricted on $T$ interpolates the values $v_{K_m}$, $v_{K_n}$ and $0$.
	\item $\overline{T}=\overline{T_{m,n}}$ for some $(m,n)\in \mathcal{N}$, and $v_{\mathcal{T}}$ restricted on $T$ interpolates the values $v_{K_m}$, 0 and~0.
\end{enumerate}
We denote by $(\lambda_k)_{1\leq k \leq 3}$ the barycentric coordinates in $T$, such that for all $x \in T$,
\[  v_{\mathcal{T}}(x) = v_{K_1} \lambda_1(x)+ v_{K_2} \lambda_2(x)+ v_{K_3} \lambda_3(x), \]
where some $v_{K_m}$ could be equal to 0 if $x_{K_m} \in \partial \Omega$. We can then write that for some $1 \leq j \leq 3$, 
\[  \nabla  v_{\mathcal{T}}(x) = \sum_{1 \leq m \neq j \leq 3} (v_{K_m}-v_{K_j}) \nabla \lambda_m(x) \]
for all $x \in T$ since $\sum_{k=1}^3 \lambda_k = 1$. For the (i) case, for $1 \leq m \neq n \leq 3$ we compute 
\begin{align*}
\int_{T_{mn}} |\nabla v_{\mathcal{T}} |^2 & = |T_{mn}| \left| \sum_{j\neq m} (v_{K_j}-v_{K_m}) \nabla \lambda_j  \right|^2 \\
& \leq 2 |T_{mn}| \sum_{j \neq m} \left(\frac{D_{\sigma_{jm}}v_{\mathcal{M}}}{d_{\sigma_{jm}}}  \right)^2 |\nabla \lambda_j|^2 d_{\sigma_{jm}}^2 \\
& \leq 2 \left( \max_{j \neq m} |\nabla\lambda_j|^2 d_{\sigma_{j m}}^2\right) \left( |T_{mn}| \left( \frac{D_{\sigma_{nm}}v_{\mathcal{M}}}{d_{\sigma_{nm}}} \right)^2 + \underbrace{|T_{mn}|}_{\lesssim |T_{m \ell}|} \left(\frac{D_{\sigma_{\ell m}}v_{\mathcal{M}}}{d_{\sigma_{\ell m}}}  \right)^2 \right)
\end{align*}
for $1 \leq \ell \leq 3$ and $\ell \neq m,n$. We either have
\[ \max_{j \neq m} |\nabla\lambda_j|^2 d_{\sigma_{j m}}^2 =  |\nabla\lambda_n|^2 d_{\sigma_{n m}}^2 = \frac{d_{\sigma_{\ell m}}^2}{4|T|^2} d_{\sigma_{n m}}^2  \]
if $j=n$, or
\[ \max_{j \neq m} |\nabla\lambda_j|^2 d_{\sigma_{j m}}^2 =  |\nabla\lambda_\ell|^2 d_{\sigma_{\ell m}}^2 = \frac{d_{\sigma_{n m}}^2}{4|T|^2} d_{\sigma_{\ell m}}^2  \]
if $j=\ell$, hence in both cases
\[   \max_{j \neq m} |\nabla\lambda_j|^2 d_{\sigma_{j m}}^2 = \frac{d_{\sigma_{n m}}^2 d_{\sigma_{\ell m}}^2}{4|T|^2} \lesssim 1    \]
by mesh assumptions. We treat the (ii) and (iii) cases similarly: for $j=\ell\neq m \neq n$ so that $x_{K_j}=x_T \in \partial \Omega$, and we let $D_{\sigma_{jm}}v_{\mathcal{M}} \coloneqq|v_{K_m}|$ and $d_{\sigma_{jm}} \coloneqq d_{\sigma_{nm}}$. The conclusion follows, in view of the definition of $|\cdot|_{1,2,\mathcal{M}}$, by summing over the $T_{mn}$ composing $T$ for $1\leq m \neq n\leq 3$, and invoking that
\[ |T_{mn}|=\frac{d_{K_m,\sigma_{mn}} |\sigma^\top_{mn}|}{2},\]
then by summing over $T\in \mathcal{T}$.
\end{proof}

\begin{lemma} \label{lemma_eq_discrete_Lp}
Let $p>2$, then there exists a constant $C>0$ independent of $p$ such that
\[ \| v_{\mathcal{M}}  \|_{0,p} \leq C \| v_{\mathcal{T}} \|_{L^p(\Omega)}.  \]
\end{lemma}
\begin{proof}
Let $T\in\mathcal{T}$. By Hölder's inequality, we know that
\[ \int_T |v_{\mathcal{T}}|^2 \leq \left( \int_T |v_{\mathcal{T}}|^p  \right)^{\frac{2}{p}}  |T|^{1-\frac{2}{p}}. \]
On the other hand, keeping in mind that some values of $v_{K_m}$ for $1 \leq m \leq 3$ might be equal to 0 (as in the proof of Lemma \ref{lemma_eq_discrete_H1_dual}) and denoting $(\lambda_k)_{1\leq k \leq 3}$ the barycentric coordinates in $T$, we can write that
\[ \int_T |v_{\mathcal{T}}|^2 = V_T^\top \mathbb{M}_T V_T, \quad \text{with} \ V_T \coloneqq \left( \begin{array}{c} v_{K_1} \\ v_{K_2} \\ v_{K_3} \end{array} \right) \ \text{and} \ \left[ \mathbb{M}_T \right]_{mn} \coloneqq \lambda_m \lambda_n  \]
on $T$. As the smallest eigenvalue of $ \mathbb{M}_T$ is equal to $|T|/12$, we have
\[ \int_T |v_{\mathcal{T}}|^2  \geq  \frac{|T|}{12} \left(v_{K_1}^2+v_{K_2}^2+v_{K_3}^2 \right).  \]
We now let $\mathcal{N}_T \subset \mathcal{N}$ be the set such that $\overline{T}=\cup_{(m,n)\in \mathcal{N}_T} \overline{T_{mn}}$, hence we can write that
\[  \int_T |v_{\mathcal{M}}|^2 = \sum_{(m,n)\in\mathcal{N}_T} |T_{mn}| |v_{K_m}|^2.  \]
Since $|T|\geq |T_{mn}|$ for any $(m,n)\in\mathcal{N}_T$, we infer that
\[ \frac{1}{24} \sum_{(m,n)\in\mathcal{N}_T} |T_{mn}| |v_{K_m}|^2 \leq \int_T |v_{\mathcal{T}}|^2. \]
Gathering these inequalities, we compute
\begin{align*}
|T|^{\frac{p}{2}-1} \int_T |v_{\mathcal{T}}|^p & \geq \left( \int_T |v_{\mathcal{T}}|^2 \right)^{\frac{p}{2}}  \geq \left( \frac{1}{24} \sum_{(m,n)\in\mathcal{N}_T} |T_{mn}| |v_{K_m}|^2 \right)^{\frac{p}{2}} \\
& \geq \frac{1}{24^{\frac{p}{2}}} \sum_{(m,n)\in\mathcal{N}_T} |T_{mn}|^{\frac{p}{2}} |v_{K_m}|^{p} =  \frac{1}{24^{\frac{p}{2}}} \sum_{(m,n)\in\mathcal{N}_T} |T_{mn}|^{\frac{p}{2}-1} \int_{T_{mn}} |v_{K_m}|^{p} \\
& \geq \frac{1}{24^{\frac{p}{2}}} (C')^{\frac{p}{2}-1} |T|^{\frac{p}{2}-1} \int_{T} |v_{\mathcal{M}}|^{p}
\end{align*}
where we have used that from the mesh assumptions $|T_{mn}| \geq C'|T|$ for a constant $C'\in \left(0,1\right]$ and for any $(m,n)\in\mathcal{N}$. Finally, summing over $T \in \mathcal{T}$, we get that
\[ \| v_{\mathcal{M}} \|_{0,p} \leq  2\sqrt{6} (C')^{\frac{1}{p}-\frac{1}{2}}  \| v_{\mathcal{T}} \|_{L^p(\Omega)}  \leq   2\sqrt{6} (C')^{-\frac{1}{2}} \| v_{\mathcal{T}} \|_{L^p(\Omega)}\]
which gives the result.
\end{proof}

\begin{corollary} \label{corollary_discrete_Sobolev}
There exists a constant $C>0$ independent of $p>2$ such that
\[ \| v_{\mathcal{M}}  \|_{0,p} \leq C \sqrt{p} |v_{\mathcal{M}} |_{1,2,\mathcal{M}}.  \]
\end{corollary}
\begin{proof}
This property is a direct consequence of Lemma \ref{lemma_eq_discrete_H1_dual} and Lemma \ref{lemma_eq_discrete_Lp}, as well as the well known Sobolev inequality on convex domains (see for instance \cite[Theorem 4.2]{Mizuguchi2017}) which states that for all~$\varphi \in H^1(\Omega)$, 
\[ \| \varphi \|_{L^p(\Omega)} \leq C_{\Omega} \sqrt{p}    \| \varphi \|_{H^1(\Omega)}. \]
Applying this inequality with $\varphi=v_{\mathcal{T}} \in H^1_0(\Omega)$ alongside Poincaré inequality gives the result.
\end{proof}

From classical Sobolev embeddings, we know that $\| f \|_{L^{\infty}(\Omega)} \leq C_{\Omega,\eps} \| f \|_{H^{1+\eps}(\Omega)}$ in dimension $d=2$ for any $\eps>0$, however such estimate do not hold for $\eps=0$. However, taking advantage of the discrete structures of our underlying spaces, in particular relying on inverse estimate for discrete Lebesgue space as well as discrete Sobolev embeddings from Corollary \ref{corollary_discrete_Sobolev}, one can show non-uniform discrete Sobolev embeddings for $\| \cdot \|_{0,\infty}$ norm.
 
\begin{corollary} \label{discrete_infty_lemma}
There exists a constant $C=C_{\Omega,\Theta_{\mathcal{V}}}>0$ such that 
\[ \| v_{\mathcal{M}} \|_{0,\infty} \leq C_{\Omega,\Theta_{\mathcal{V}}}  |\ln h_{\mathcal{M}} |^{\frac12} | v_{\mathcal{M}} |_{1,2,\mathcal{M}}   \]
\end{corollary}
\begin{proof}
There exists $K_0 \in \mathcal{T}$ such that $\| v_{\mathcal{M}} \|_{0,\infty}= |v_{K_0}|$, so
\[ \| v_{\mathcal{M}} \|_{0,\infty} \leq C h_{\mathcal{M}}^{-\frac{2}{p}} \| v_{\mathcal{M}} \|_{0,p}   \]
where $C=C(\Theta_{\mathcal{V}})$ depends only on the regularity of the mesh. From Corollary \ref{corollary_discrete_Sobolev}, we then get that
\[ \| v_{\mathcal{M}} \|_{0,\infty} \leq C h_{\mathcal{M}}^{-\frac{2}{p}} \sqrt{p} |v_{\mathcal{M}} |_{1,2,\mathcal{M}}  \]
so taking $p=|\ln h_{\mathcal{M}}|$ gives the result.
\end{proof}

We also state the discrete equivalent of the usual integration by parts formula in the context of finite-volume schemes: for $v_{\mathcal{M}}$, $\widetilde{v}_{\mathcal{M}} \in X_{\mathcal{V},0}$ such that $v_K=\widetilde{v}_K=0$ for $K \in \mathcal{M}_{\interior}$, we have
\[ \sum_{K \in \mathcal{M}}\sum_{\sigma \in \mathcal{E}_K \cap \mathcal{E}_{\interior}} \frac{|\sigma|}{d_{K|L}} (v_K-v_L)\widetilde{v}_K  =  \sum_{\sigma \in \mathcal{E}_{\interior}} \frac{|\sigma|}{d_{K|L}} (v_K-v_L)(\widetilde{v}_K - \widetilde{v}_L) \]
rearranging the sum on the left-hand side.

\subsection{Interpolants}

In this section we compare the pointwise interpolant $P_h$ with the mean interpolant $\pi_{\mathcal{M}}: L^1(\Omega) \rightarrow X_{\mathcal{V},0}$ defined by
\[ \left. (\pi_{\mathcal{M}}\varphi) \right|_K \coloneqq \frac{1}{|K|} \int_K \varphi(x) \dd x  \]
for all $K \in \mathcal{M}$ and $\varphi \in L^1(\Omega)$.

\begin{lemma} \label{lemma_interpolant_diff}
Let $\varphi \in H^3\cap H^1_0(\Omega)$, then
\[ \| P_{\mathcal{M}} \varphi - \pi_{\mathcal{M}} \varphi \|_{0,2} \leq C h_{\mathcal{M}} \| \varphi \|_{H^3(\Omega)}.  \]
\end{lemma}
\begin{proof}
Let $K \in \mathcal{M}$. This property is a direct consequence of the Taylor expansion
\[  \varphi(x)-\varphi(x_K)=\int_0^1 \nabla \varphi (\theta x_K +(1-\theta)x) \cdot (x-x_K) \dd \theta  \]
for any $x \in K$. Integrating over $K$ and dividing by $|K|$, we get that
\begin{align*}
\left| \varphi(x_K)-\frac{1}{|K|}\int_K \varphi(y) \dd y \right| & = \frac{1}{|K|}  \left| \int_K \int_0^1 \nabla \varphi (\theta x_K +(1-\theta)x) \cdot (x-x_K) \dd \theta  \dd y \right| \\
& \lesssim h_K \| \nabla \varphi \|_{L^{\infty}(\Omega)} \lesssim h_{\mathcal{M}} \| \varphi \|_{H^3(\Omega)}
\end{align*}   
by usual Sobolev embeddings. Taking the square, multiplying by $|K|$, summing over $K \in \mathcal{M}$ and taking the square root gives the result.
\end{proof}

\subsection{Conserved quantities of the finite volume scheme} \label{section_conserved_quantities}
It is useful to remark that multiplying equation \eqref{FVGP} taken on a control volume $K$ by $|K|\overline{u_K}$, summing over $K \in \mathcal{M}$ and taking the imaginary part, one infer that
\begin{align*}
\frac12 \frac{\dd }{\dd t} \| u_{\mathcal{M}}(t)\|_{0,2}^2 & = \Im \sum_{K\in\mathcal{T}} |K| \left(  \overline{u_K(t)}   (A_{\mathcal{V}}u_{\mathcal{M}}(t))_K + \gamma |u_K(t)|^4 + V(t,x_K) |u_K(t)|^2\right) \\
& =\Im \sum_{K\in\mathcal{T}} \left( \sum_{\sigma \in \mathcal{E}_{\interior}\cap \mathcal{E}_K} \frac{|\sigma|}{d_{K|L}} |u_K(t)-u_L(t)|^2 \right)= 0
\end{align*}   
from discrete integration by parts, which leads to the conservation of the discrete $L^2$-norm 
\[ \| u_{\mathcal{M}}(t)\|_{0,2} = \| u_{\mathcal{M}}(0)\|_{0,2} \]
 for all $t \in \left[0,T \right]$. The same way, multiplying  equation \eqref{FVGP} by $|K|\partial_t \overline{u_K}$, summing over $K \in \mathcal{M}$ and taking this time the real part, we end up with the discrete energy balance law
\[ \frac{\dd }{\dd t} H_{\mathcal{V}}(t) = \sum_{K\in \mathcal{M}} |K| \partial_t V(t,x_K) |u_K(t)|^2    \]
where
\[ H_{\mathcal{V}}(t) \coloneqq   | u_{\mathcal{M}}(t) |_{1,2,\mathcal{M}}^2 + \frac{\gamma}{2} \| u_{\mathcal{M}}(t) \|_{0,4}^4 +  \|(P_{\mathcal{M}}V)(t) u_{\mathcal{M}}(t)\|_{0,2}^2. \]
Integrating this equation over time on $\left[0,T\right]$, from the assumption that $V$ and $\partial_t V$ are bounded by below and by Cauchy-Schwarz inequality (using the discrete mass conservation property) we infer the uniform (with respect to $h_{\mathcal{M}}$) bound
\begin{equation} \label{eq_uniform bound_U}
\sup_{t\in \left[0,T \right]} | u_{\mathcal{M}}(t) |_{1,2,\mathcal{M}} \leq C.
\end{equation}
for our numerical scheme, which will be useful in the upcoming analysis.

\section{Error estimates} \label{section_error}
We turn to the proof of Theorem \ref{theorem_convergence}. We denote the error by $e_{\mathcal{M}}=(e_K)_{K \in \mathcal{M}}$, defined by
\[ e_K(t) \coloneqq u(t,x_K)-u_K(t)  \]
for all $t\in \left[0,T\right]$ and any $K\in\mathcal{M}$. Integrating equation \eqref{gross_pitaevskii} over a control volume~$K$ and using Green formula on the Laplace operator we end up with
\[ i\int_K \partial_t u(t,x) \dd x + \sum_{\sigma \in \mathcal{E}} \int_{\sigma} \nabla u(t,x) \cdot \nu_{K,\sigma} \dd S = \gamma \int_K  |u(t,x)|^2 u(t,x) \dd x + \int_{K} V(t,x) u(t,x) \dd x \]
where $ \nu_{K,\sigma}$ denotes the outward-pointing unit normal vector to $\sigma$. Dividing this equation by $|K|$ and substracting equation \eqref{FVGP}, we infer the following equation for the error term
\begin{multline} \label{eq_error_term_equation}
i e_K'(t) + (A_{\mathcal{V}} e_{\mathcal{M}})_K(t)= \gamma |u(t,x_K)|^2 e_K(t) +\gamma u(t,x_K) u_K(t) \overline{e_K(t)} + \gamma |u_K(t)|^2 e_K(t)  \\
+ V(t,x_K)e_K(t) - i \mathcal{S}_K(t) - \mathcal{R}_K(t) +\gamma \mathcal{U}_K(t)    + \mathcal{W}_K(t)  ,
\end{multline}
where
\[ \left\{ 
\begin{aligned}
&\mathcal{S}_K(t)\coloneqq \frac{1}{|K|} \int_K \partial_t u(t,x) \dd x - \partial_t u(t,x_K),\\
&\mathcal{R}_K(t)\coloneqq \frac{1}{|K|} \sum_{\sigma \in \mathcal{E}_{\interior}} \left( \int_{\sigma} \nabla u(t,x) \cdot \nu_{K,\sigma} \dd S - \frac{|\sigma|}{d_{K|L}} \left(u(t,x_L)-u(t,x_K)\right) \right) \\
& \quad \quad \quad \quad + \frac{1}{|K|} \sum_{\sigma \in \mathcal{E}_{\exterior}}  \int_{\sigma} \nabla u(t,x) \cdot \nu_{K,\sigma} \dd S ,\\
&\mathcal{U}_K(t)\coloneqq \frac{1}{|K|} \int_K  |u(t,x)|^2 u(t,x) \dd x - |u(t,x_K)|^2 u(t,x_K), \\
& \mathcal{W}_K(t)\coloneqq \frac{1}{|K|} \int_K V(t,x) u(t,x) \dd x - V(t,x_K) u(t,x_K).
\end{aligned} \right. \]
Note that the error terms $\mathcal{S}$, $\mathcal{R}$, $\mathcal{U}$ and $\mathcal{W}$ only depends on the continuous solution $u$. Similarly to the computations performed in Section \ref{section_conserved_quantities}, multiplying equation \eqref{eq_error_term_equation} by $|K| \overline{e_K(t)}$, summing over~$K \in \mathcal{M}$ and taking the imaginary part, we end up, from several simplifications from real terms (in particular the discrete Laplace operator term which simplifies from discrete ingegration by parts), with
\begin{align*} \label{eq_error_integrated}
\frac12 \frac{\dd}{\dd t} \| e_{\mathcal{M}}(t)\|_{0,2}^2 & =  \Im \sum_{K\in \mathcal{T}} |K| \overline{e_K(t)} \left( \gamma u(t,x_K) u_K(t) \overline{e_K(t)} - i \mathcal{S}_K(t) - \mathcal{R}_K(t) +\gamma \mathcal{U}_K(t)    + \mathcal{W}_K(t)  \right)   \\
& \leq \gamma \| e_{\mathcal{M}}(t)\|_{0,2}^2 \|u(t)\|_{L^{\infty}(\Omega)} \|u_{\mathcal{M}}(t)\|_{0,\infty} + \| e_{\mathcal{M}}(t)\|_{0,2}  \| \mathcal{S}_{\mathcal{M}}(t)\|_{0,2}  \\
& \quad - \Im \sum_{K\in \mathcal{M}} |K| \overline{e_K(t)} \mathcal{R}_K(t) + \gamma \| e_{\mathcal{M}}(t)\|_{0,2} \| \mathcal{U}_{\mathcal{M}}(t)\|_{0,2} + \| e_{\mathcal{M}}(t)\|_{0,2} \| \mathcal{W}_{\mathcal{M}}(t)\|_{0,2} \\
& \eqqcolon \gamma I_1(t) + I_2(t) - I_3(t) + \gamma I_4(t) + I_5(t)
\end{align*}   
from Cauchy-Schwarz and Hölder inequalities, so we then need to estimate each of the terms on the right hand side. We start with $I_1$, using first the assumptions on $u$ which provides that~$\|u(t)\|_{L^{\infty}(\Omega)} \lesssim \| u \|_{L^{\infty}_T H^2(\Omega)}$ from Sobolev embeddings. We also use the discrete Sobolev embedding from Lemma~\ref{discrete_infty_lemma} alongside the $|\cdot|_{1,2,\mathcal{M}}$-uniform bound from the energy estimate, namely equation \eqref{eq_uniform bound_U}, to get that 
\[ \|u_{\mathcal{M}}(t)\|_{0,\infty} \lesssim |\ln h_{\mathcal{M}} |^{\frac12} \sup_{t \in \left[0,T \right]} |u_{\mathcal{M}}|_{1,2,\mathcal{M}} \lesssim |\ln h_{\mathcal{M}} |^{\frac12},\]
for all $t\in \left[0,T\right]$, hence
\[    I_1(t) \lesssim |\ln h_{\mathcal{M}}|^{\frac12} \| e_{\mathcal{M}}(t)\|_{0,2}^2.  \]
To get a uniform bound on $I_2$, we simply use Lemma \ref{lemma_interpolant_diff} with $\varphi=\partial_t u$, so that
\[ \| \mathcal{S}_{\mathcal{M}}(t)\|_{0,2} \lesssim h_{\mathcal{M}} \| \partial_t u \|_{L^{\infty}_T H^3(\Omega)}.  \]
We then infer that $I_2(t) \lesssim h_{\mathcal{M}} \| e_{\mathcal{M}}(t)\|_{0,2}$ from the assumptions on $u$ and $V$, as from the expression of equation \eqref{gross_pitaevskii}, taking the $H^3$-norm we have
\begin{align*}
\| \partial_t u \|_{ L^{\infty}_T H^3(\Omega)} & \leq \| \Delta u \|_{L^{\infty}_T H^3(\Omega)} +   \| \gamma |u|^2 u \|_{L^{\infty}_T H^3(\Omega)} + \| V u \|_{L^{\infty}_T H^3(\Omega)} \\
& \lesssim \| u \|_{L^{\infty}_T H^5(\Omega)} + \|u \|_{L^{\infty}_T H^3(\Omega)}^3 + \| V \|_{L^{\infty}_T H^3(\Omega)} \| u \|_{L^{\infty}_T H^3(\Omega)} \lesssim C
\end{align*}  
from Sobolev algebra rule and Sobolev embeddings. To estimate $I_3$, we remark that as $u(t,x_K)=u_K=0$ hence $e_K(t)=0$ if $K \in \mathcal{M}_{\exterior}$, we can restrict our attention to $\sigma \in \mathcal{E}_{\interior}$, and we denote 
\[ R_{\sigma}(t)\coloneqq \frac{1}{|\sigma|}\int_{\sigma} \nabla u(t,x) \cdot \nu_{K,\sigma} \dd S - \frac{1}{d_{K|L}} \left(u(t,x_L)-u(t,x_K)\right) \]
for any $\sigma=K|L \in \mathcal{E}_{K}\cap \mathcal{E}_{\interior}$ and $K$, $L \in \mathcal{M}$, so that $\mathcal{R}_K(t)= \frac{1}{|K|} \sum_{\sigma \in \mathcal{E}_{\interior}} |\sigma| R_{\sigma}(t) $. We then perform a Taylor expansion in space, writing that
\[u(x_L) - u(x_K)=(x_L-x_K) \cdot \nabla u(x_K) + \frac12 \int_0^1 H(u)(\theta x_K+(1-\theta)x_L)(x_L-x_K)\cdot(x_L-x_K) \dd \theta.    \]
On the other hand, from the mean value theorem we know that there exists $\eta_{\sigma} \in \sigma$
\[ \frac{1}{|\sigma|} \int_{\sigma} \nabla u(t,x) \cdot \nu_{K,\sigma} \dd S = \nabla u(t,\eta_{\sigma}) \cdot \nu_{K,\sigma} \]
thanks to the regularity of $u$, so we can also perform a Taylor expansion
\[ \nabla u(t,\eta_{\sigma}) \cdot \nu_{K,\sigma}=\nabla u (t,x_K) \cdot \nu_{K,\sigma} + \int_0^1 \nabla u (t,\theta x_K +(1-\theta)\eta_{\sigma})\cdot(\eta_{\sigma}-x_K) \dd \theta.  \]
Substracting both Taylor expansions, from simplifications as $x_L-x_K=\nu_{K,\sigma} d_{K|L}$, we compute
\begin{align*}
|R_{\sigma}| & \leq \left| \int_0^1 \nabla u (t,\theta x_K +(1-\theta)\eta_{\sigma})\cdot(\eta_{\sigma}-x_K) \dd \theta \right. \\
& \quad \quad \left. + \frac{1}{2 d_{K|L}} \int_0^1 H(u)(\theta x_K+(1-\theta)x_L)(x_L-x_K)\cdot(x_L-x_K) \dd \theta  \right| \\
& \leq\| \nabla u(t) \|_{L^{\infty}(\Omega)} |\eta_{\sigma}-x_K| + \| H(u(t))\|_{L^{\infty}(\Omega)} \frac{|x_L-x_K|^2}{2  d_{K|L}} \\
& \lesssim h_{\mathcal{M}} \|u\|_{L^{\infty}_T H^4(\Omega)}.
\end{align*}
We can then write, rearranging the sum and using that $R_{K,\sigma}(t)=-R_{L,\sigma}(t)$ so $|R_{\sigma}| = |R_{K,\sigma}(t)|=|R_{L,\sigma}(t)|$ for any edge $\sigma=K|L \in \mathcal{E}_{\interior}$, that 
\begin{align*}
|I_3(t)| &= | \sum_{\sigma \in \mathcal{E}_{\interior}} |\sigma| (e_K(t)R_{K,\sigma}(t) + e_L(t) R_{L,\sigma}(t)| \\
			& \leq \sum_{\sigma \in \mathcal{E}_{\interior}} |\sigma| D_{\sigma} e_{\mathcal{M}}(t) |R_{\sigma}(t)| \\
			& \leq \left( \sum_{\sigma \in \mathcal{E}_{\interior}} \frac{|\sigma|}{d_{K|L}} |D_{\sigma}e_{\mathcal{M}}(t)|^2   \right)^{\frac12} \left( \sum_{\sigma \in \mathcal{E}_{\interior}} |\sigma| d_{K|L} |R_{\sigma}(t)|^2  \right)^{\frac12} \\
			& \lesssim |e_{\mathcal{M}}|_{1,2,\mathcal{M}} \sqrt{|\Omega|} h_{\mathcal{M}}  \\
			& \lesssim \left( \|u\|_{L^{\infty}_T H^4(\Omega)} + \sup_{t\in\left[0,T\right]} |u_{\mathcal{M}}|_{1,2,\mathcal{M}} \right) h_{\mathcal{M}} \lesssim h_{\mathcal{M}} 
\end{align*}
using the discrete bound \eqref{eq_uniform bound_U}. Finally for $I_4$, we directly use Lemma \ref{lemma_interpolant_diff} with $\varphi=|u|^2 u$ which leads with Sobolev product rule to the bound
\[ I_4(t) \lesssim h_{\mathcal{M}} \|u\|_{L^{\infty}_T H^3(\Omega)}^3 \| e_{\mathcal{M}}(t)\|_{0,2}, \]
and we perform similarly for $I_5$ using Lemma \ref{lemma_interpolant_diff} with $\varphi=V u $, from which we infer
\[ I_5(t) \lesssim h_{\mathcal{M}} \|V\|_{L^{\infty}_T H^3(\Omega)} \|u\|_{L^{\infty}_T H^2(\Omega)} \| e_{\mathcal{M}}(t)\|_{0,2}. \]
Gathering the bounds on $I_1$, $I_2$, $I_3$, $I_4$ and $I_5$, we then get that
\[  \frac12 \frac{\dd}{\dd t} \| e_{\mathcal{M}}(t)\|_{0,2}^2 \lesssim |\ln h_{\mathcal{M}} |^{\frac12}  \| e_{\mathcal{M}}(t)\|_{0,2}^2 + h_{\mathcal{M}} \| e_{\mathcal{M}}(t)\|_{0,2} + h_{\mathcal{M}}, \]
so from Gronwall lemma we infer that
\[ \| e_{\mathcal{M}}(t)\|_{0,2} \lesssim h_{\mathcal{M}}^{\frac12} e^{C(T)  |\ln h_{\mathcal{M}} |^{\frac12}}  \]
for a constant $C(T)>0$. As the function $h_{\mathcal{M}}^{\eps} e^{C(T)  |\ln h_{\mathcal{M}} |^{\frac12}} \underset{h_{\mathcal{M}}\to 0}{\longrightarrow}0$ for any $\eps>0$, we get the result.

\section{Numerics} \label{section_numerics}

In this section, we report numerical results of the proposed numerical method. We take $\Omega$ as a polygonal approximation of the open disk $\mathcal{D}(0,2)$ of center $(0,0)\in \R^2$ and of radius~2. To generate an admissible Delaunay triangulation of $\Omega$ in the sens of \cite{eymard2000}, we use the \textsc{GMSH} software \cite{gmsh} with the "Frontal-Delaunay" option. We can then build the dual Voronoï mesh used for the numerical scheme \eqref{FVGP}, see Figure \ref{fig:mesh}.

\begin{figure}[h]
	\centering
		\includegraphics[width=1\textwidth,trim = 25cm 25cm 25cm 25cm, clip]{mesh.png}	
	\caption{Voronoï meshes (black) of $\Omega$ and their dual Delaunay triangulations (gray), with 71 Voronoï cells (117 triangles) on the left and 222 Voronoï cells (397 triangles) on the right.}
	\label{fig:mesh}
\end{figure} 

To carry out our simulations, we also need to perform a time discretization. We employ a Strang splitting integration, at it is fairly easy to implement and robust, and naturally preserves the mass of the numerical solution. Fixing a number of iterations $N \in \N^*$ and taking a time step-size $\tau=T/N$, it ends up composing the linear flow 
\[ \Phi_A^{\tau} (v_{\mathcal{M}}) \coloneqq e^{i\tau A_{\mathcal{V}}} v_{\mathcal{M}}  \]
and the nonlinear flow
\[ \Phi_B^{\tau} (w_{\mathcal{M}})\coloneqq e^{-i\tau \gamma |w_{\mathcal{M}}|^2 } e^{-i\tau \int_0^{\tau} (P_{\mathcal{M}}V)(s) \dd s} w_{\mathcal{M}}. \]
 Denoting one step of time integration as
\[  \Phi_S^{\tau}(z_{\mathcal{M}})\coloneqq \Phi_B^{\frac{\tau}{2}} \circ \Phi_A^{\tau} \circ \Phi_B^{\frac{\tau}{2}}(z_{\mathcal{M}}), \]
this leads to the approximation for the final time
\[ u_{\mathcal{M}}(T) \simeq \underbrace{\Phi_S^{\tau} \circ \ldots \circ \Phi_S^{\tau}}_{N \ \text{times}}( P_{\mathcal{M}} u_0).  \]
Note that for numerical efficiency, the linear flow $\Phi_A^{\tau}$ is computed using a two-order Padé approximant
\[ e^{i \tau A_{\mathcal{V}}} \simeq \left( \Id_{d_{\mathcal{M}}} - i \frac{\tau}{2}A_{\mathcal{V}}  \right)^{-1} \left( \Id_{d_{\mathcal{M}}} + i \frac{\tau}{2}A_{\mathcal{V}}  \right)  \]
which also preserves the $L^2(\Omega)$-norm of the numerical solution. We emphasize that the matrix~$A_{\mathcal{V}}$ is a sparse matrix, so that the resolution of the linear system $A_{\mathcal{V}}X=B$ can be efficiently precomputed using its sparse LU decomposition. Note that all the codes used in this section are available on the \textsc{Gitlab} page \url{https://plmlab.math.cnrs.fr/chauleur/codes/-/tree/main/FVGP_codes}.

\subsection{Accuracy test}
We first restrict our attention to the linear case, taking $\gamma=0$. In this case, for a radial harmonic trap $V(r)=\omega^2 |x|^2$ and initial condition taken as a radial centered Gaussian function~$u_0(x)=e^{-4|x|^2}$, equation \eqref{gross_pitaevskii} has an explicit Gaussian solution of the form
\[ u(t,x)= \exp \left( -\alpha(t)|x|^2 - 4i\int_0^t \alpha(s) \dd s \right), \quad \alpha(t)=-\frac{i}{4} \frac{16i \cos(2\omega t)-2\omega \sin(2\omega t)}{\frac{8 i}{\omega} \sin(2\omega t) + \cos( 2\omega t)},   \]
for all $t\in \R$ and $x \in \R^2$. We then turn to the nonlinear case taking $\gamma=1$, with the same radial harmonic trap $V$. To get a reference solution, we look for the ground state solution of the form $u(t,x)=e^{-i\mu t} \phi(x)$, where $\phi$ is the unique solution of the nonlinear equation
\[ \mu \phi= -\Delta \phi + \gamma |\phi|^2\phi + V \phi  \] 
with $\| \phi \|_{L^2(\Omega)}=1$, and the chemical potential $\mu$ is given by
\[ \mu = \int_{\Omega} |\nabla \phi |^2 +  \gamma \int_{\Omega} | \phi |^4 + \int_{\Omega} V | \phi |^2. \]
To compute the initial state $u_0(x)=\phi(x)$, we perform a normalized gradient flow in 1D radial coordinates like in \cite{bao2004} with a very precise stepsize, that we then interpolate on our Voronoï grid.

\begin{figure}[h]
	\centering
		\includegraphics[width=1\textwidth,trim = 10cm 10cm 10cm 10cm, clip]{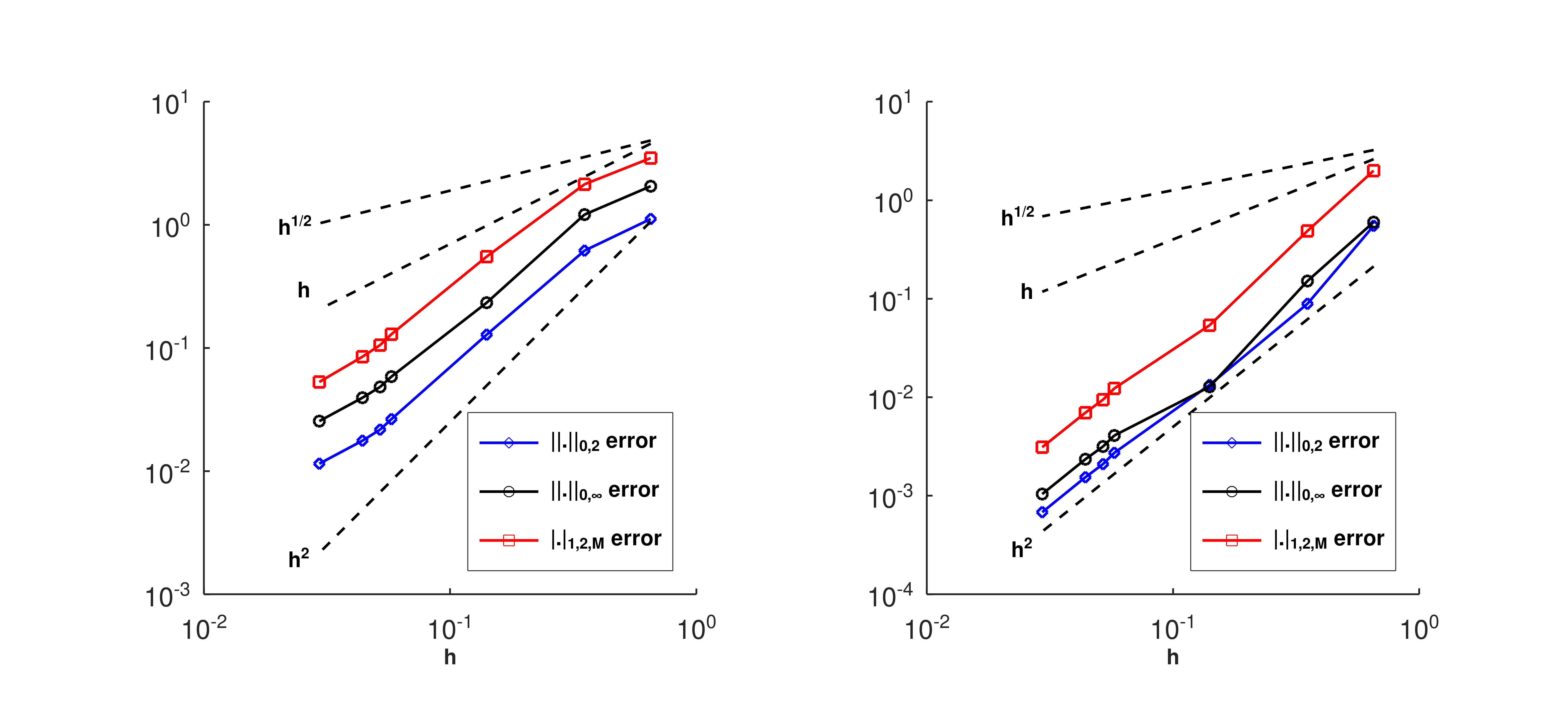}	
	\caption{Convergence of the numerical errors in $\| \cdot \|_{0,2}$, $\| \cdot \|_{0,\infty}$ and $|\cdot|_{1,2,\mathcal{M}}$ norms for the linear case \textit{(left)} and the nonlinear ground state \textit{(right)} as~$h_{\mathcal{M}} \rightarrow 0$.}
	\label{fig:error}
\end{figure}

For both cases, we compare our numerical result with the true solution at $T=1$ taking $N=2000$, and generating Voronoï meshes with several diameters $h_{\mathcal{M}}$ between $0.5$ and $0.05$ (which results from Delaunay triangulations with between approximately 100 triangles and 60000 triangles). We then compute numerical errors of the form
\[ \err(h_{\mathcal{M}}) \coloneqq \| u(T) - u_{\mathcal{M}}(T)  \|   \]
 in respectively $\| \cdot \|_{0,2}$, $\| \cdot \|_{0,\infty}$ and $|\cdot|_{1,2,\mathcal{M}}$ norms, for both the linear and nonlinear cases, and plot the results in Figure~\ref{fig:error}. The numerical results show an in-between first-order and second-order convergence in all three norms, which is better than the rate proven theoretically in Theorem \ref{theorem_convergence}. 

\subsection{Vortex nucleation}
\begin{figure}[h]
	\centering
		\includegraphics[width=1\textwidth,trim = 40cm 40cm 40cm 40cm, clip]{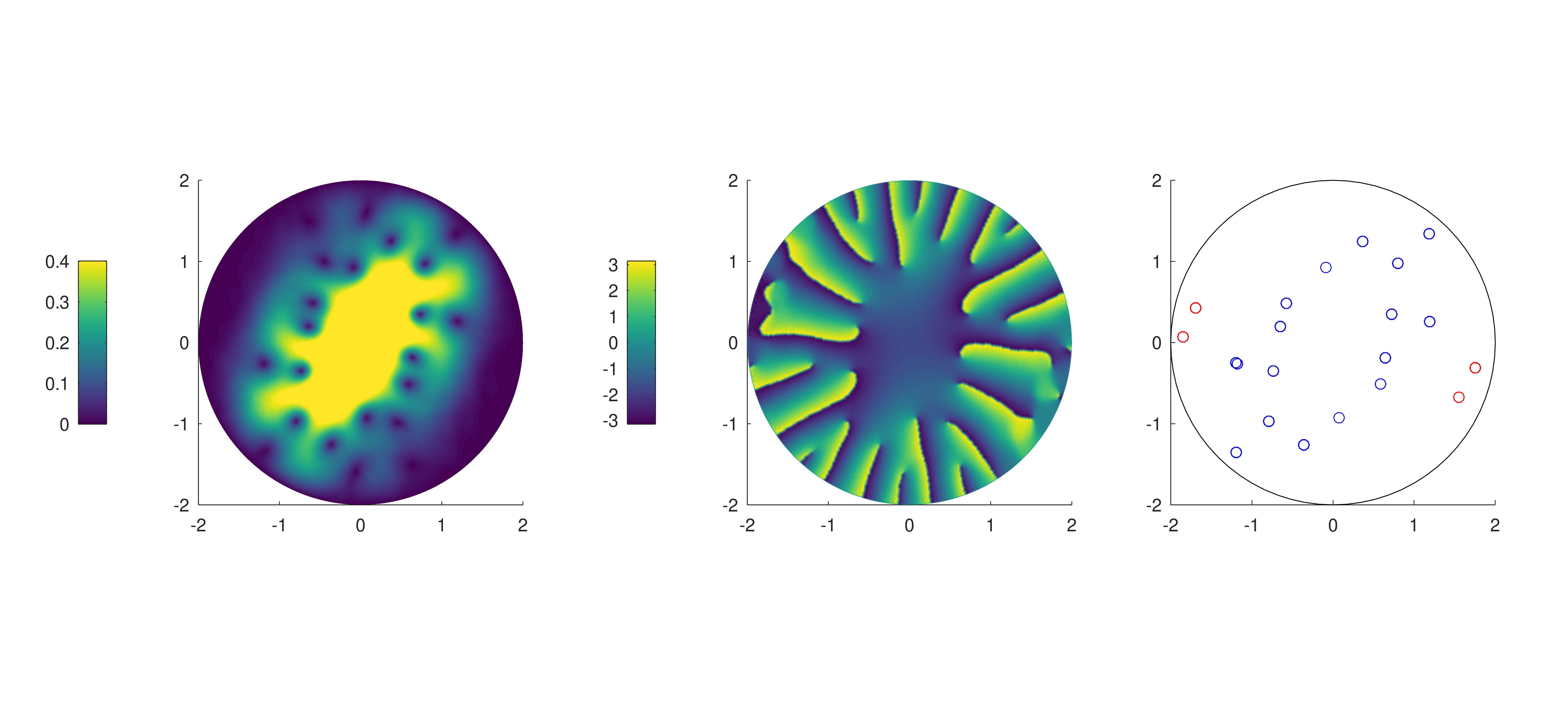}	
	\caption{Plots of the density $|u_{\mathcal{M}}(T)|$ and the phase $\arg(u_{\mathcal{M}}(T)) \in \left[-\pi,\pi \right[$.}
	\label{fig:vortices}
\end{figure}

To place our work within its relevant physical framework, namely, the theory of quantum turbulence, we present an application in this setting. Specifically, we consider the case where the harmonic potential is perturbed by a rotating sinusoidal stirrer of size $\epsilon$ and frequency $\omega$, which writes
\[ V(t,x)=\omega^2 r^2 \left(1+ \epsilon \cos(2\theta - \beta t) \right), \quad x=(r,\theta), \]
in polar coordinates. In the upcoming simulations we will take $\omega=10$, $\epsilon=0.2$ and $\beta=30$. We also take the nonlinearity constant $\gamma=100$ to enhance vortex nucleation. As an initial state $u_0$, we start from the ground state solution associated to equation \eqref{gross_pitaevskii} with $\epsilon=0$, which is computed performing a normalized gradient flow (similarly to the accuracy test of the previous section, see also \cite[Section 4]{Chauleur2024}). We perform our simulation with $T=5$ and time step size~$\tau=0.001$, and plot the results obtained in Figure \ref{fig:vortices}. We well observe vortex nucleation, in a similar way as in the work~\cite{lundh2003}. Note that more numerical experiments with a similar finite volume scheme alongside detailed physical motivations are provided in another paper of the author~\cite{Chauleur2024}.

\subsection*{Acknowledgements}
The author is supported by the Labex CEMPI (ANR-11-LABX-0007-01). The author is grateful to Claire Chainais-Hillairet, Clément Cancès, Guillaume Dujardin, Guillaume Ferrière, Simon Lemaire and Julien Moatti for helpful discussions about this work, in particular to Simon Lemaire for pointing out the proof of Lemma \ref{lemma_eq_discrete_H1_dual} and Lemma \ref{lemma_eq_discrete_Lp}.

\bibliographystyle{siam}
\bibliography{biblio}

\begin{thebibliography}{10}

\bibitem{Boyer2006}
{\sc B.~Andreianov, F.~Boyer, and F.~Hubert}, {\em On the finite-volume
  approximation of regular solutions of the {$p$}-{L}aplacian}, IMA J. Numer.
  Anal., 26 (2006), pp.~472--502.

\bibitem{duboscq2015stationary}
{\sc X.~Antoine and R.~Duboscq}, {\em G{P}{E}{L}ab, a {M}atlab toolbox to solve
  {G}ross–{P}itaevskii equations {I}: {C}omputation of stationary solutions},
  Computer Physics Communications, 185 (2014), pp.~2969--2991.

\bibitem{duboscq2015dynamics}
\leavevmode\vrule height 2pt depth -1.6pt width 23pt, {\em G{P}{E}{L}ab, a
  {M}atlab toolbox to solve {G}ross–{P}itaevskii equations {I}{I}: {D}ynamics
  and stochastic simulations}, Computer Physics Communications, 193 (2015),
  pp.~95--117.

\bibitem{bao2013}
{\sc W.~Bao and Y.~Cai}, {\em Mathematical theory and numerical methods for
  {B}ose-{E}instein condensation}, Kinet. Relat. Models, 6 (2013), pp.~1--135.

\bibitem{bao2004}
{\sc W.~Bao and Q.~Du}, {\em Computing the ground state solution of
  {B}ose-{E}instein condensates by a normalized gradient flow}, SIAM J. Sci.
  Comput., 25 (2004), pp.~1674--1697.

\bibitem{bradji2015}
{\sc A.~Bradji}, {\em A theoretical analysis for a new finite volume scheme for
  a linear {S}chr\"{o}dinger evolution equation on general nonconforming
  spatial meshes}, Numer. Funct. Anal. Optim., 36 (2015), pp.~590--623.

\bibitem{carles2011}
{\sc R.~Carles}, {\em Nonlinear {S}chr\"{o}dinger equation with time dependent
  potential}, Commun. Math. Sci., 9 (2011), pp.~937--964.

\bibitem{Chauleur2024}
{\sc Q.~Chauleur, R.~Chicireanu, G.~Dujardin, J.-C. Garreau, and A.~Rançon},
  {\em Numerical study of the {G}ross-{P}itaevskii equation on a
  two-dimensional ring and vortex nucleation}, 2024.
\newblock Preprint, archived at \url{https://hal.science/hal-04545713}.

\bibitem{droniou2018}
{\sc J.~Droniou, R.~Eymard, T.~Gallou\"{e}t, C.~Guichard, and R.~Herbin}, {\em
  The gradient discretisation method}, vol.~82 of Math\'{e}matiques \&
  Applications (Berlin) [Mathematics \& Applications], Springer, Cham, 2018.

\bibitem{eymard2000}
{\sc R.~Eymard, T.~Gallou\"{e}t, and R.~Herbin}, {\em Finite volume methods},
  in Handbook of numerical analysis, {V}ol. {VII}, Handb. Numer. Anal., VII,
  North-Holland, Amsterdam, 2000, pp.~713--1020.

\bibitem{gartner2019}
{\sc K.~G\"{a}rtner and L.~Kamenski}, {\em Why do we need {V}oronoi cells and
  {D}elaunay meshes? {E}ssential properties of the {V}oronoi finite volume
  method}, Comput. Math. Math. Phys., 59 (2019), pp.~1930--1944.

\bibitem{gmsh}
{\sc C.~Geuzaine and J.-F. Remacle}, {\em Gmsh: {A} 3-{D} finite element mesh
  generator with built-in pre- and post-processing facilities}, Internat. J.
  Numer. Methods Engrg., 79 (2009), pp.~1309--1331.

\bibitem{henning2017}
{\sc P.~Henning and A.~M{\aa}lqvist}, {\em The finite element method for the
  time-dependent {G}ross-{P}itaevskii equation with angular momentum rotation},
  SIAM J. Numer. Anal., 55 (2017), pp.~923--952.

\bibitem{lundh2003}
{\sc E.~Lundh, J.-P. Martikainen, and K.-A. Suominen}, {\em Vortex nucleation
  in {B}ose-{E}instein condensates in time-dependent traps}, Phys. Rev. A, 67
  (2003), p.~063604.

\bibitem{mishev1998}
{\sc I.~D. Mishev}, {\em Finite volume methods on {V}oronoi meshes}, Numer.
  Methods Partial Differential Equations, 14 (1998), pp.~193--212.

\bibitem{Mizuguchi2017}
{\sc M.~Mizuguchi, K.~Tanaka, K.~Sekine, and S.~Oishi}, {\em Estimation of
  {S}obolev embedding constant on a domain dividable into bounded convex
  domains}, J. Inequal. Appl.,  (2017), pp.~Paper No. 299, 18.

\bibitem{danaila2016}
{\sc G.~Vergez, I.~Danaila, S.~Auliac, and F.~Hecht}, {\em A finite-element
  toolbox for the stationary {G}ross–{P}itaevskii equation with rotation},
  Computer Physics Communications, 209 (2016), pp.~144--162.

\bibitem{krstulovic2020}
{\sc A.~Villois, D.~Proment, and G.~Krstulovic}, {\em Irreversible dynamics of
  vortex reconnections in quantum fluids}, Phys. Rev. Lett., 125 (2020),
  pp.~164501, 5.

\end{thebibliography}

\end{document}